\documentclass[11pt]{amsart}

\usepackage{amsopn}
\usepackage{amsmath,amsthm,amssymb, amscd, amsfonts}

\theoremstyle{theorem}
\newtheorem{teo}{\bf Theorem}[section]
\newtheorem{cor}[teo]{\bf Corollary}

\newtheorem{prop}[teo]{\bf Proposition}
\newtheorem{lema}[teo]{\bf Lemma}

\theoremstyle{definition}

\theoremstyle{remark}
\newtheorem{nota}[teo]{\bf Remark}

\numberwithin{equation}{section}

\newcommand\ad{\operatorname{ad}}
\newcommand\Ad{\operatorname{Ad}}
\newcommand\Aff{\operatorname{Aff}}

\newcommand\diag{\operatorname{diag}}

\newcommand\Exp{\operatorname{Exp}}

\newcommand\Id{\operatorname{Id}}

\newcommand\Iso{\operatorname{Iso}}
\newcommand\Lie{\operatorname{Lie}}
\newcommand\OO{\operatorname{O}}
\newcommand\SO{\operatorname{SO}}
\newcommand\Sp{\operatorname{Sp}}
\newcommand\Spin{\operatorname{Spin}}
\newcommand\SU{\operatorname{SU}}

\newcommand\trace{\operatorname{trace}}
\newcommand\Tr{\operatorname{Tr}}

\newcommand\vol{\operatorname{vol}}

\renewcommand\gg{\mathfrak{g}}
\newcommand\gh{\mathfrak{h}}

\newcommand\gm{\mathfrak{m}}
\newcommand\bbv{\mathbb{V}}
\newcommand\bbw{\mathbb{W}}
\newcommand\bbr{\mathbb{R}}
\newcommand\cd{\mathcal{D}}
\newcommand\cc{\mathcal{C}}
\newcommand\cf{\mathcal{F}}
\newcommand\cg{\mathcal{G}}
\newcommand\ck{\mathcal{K}}

\makeatletter\renewcommand\theenumi{\@roman\c@enumi}\makeatother

\begin{document}

\title[The skew-torsion holonomy theorem]{The skew-torsion holonomy theorem and naturally reductive spaces}

\author{Carlos Olmos} \address{Facultad de Matem\'atica, Astronom\'\i a y F\'\i sica, Universidad Nacional de
C\'ordoba, Ciudad Universitaria, 5000 C\'ordoba, Argentina}  \email{olmos@mate.uncor.edu \ \ \
reggiani@mate.uncor.edu}

\author{Silvio Reggiani}




\thanks {Supported by Universidad Nacional de C\'ordoba and CONICET. Partially supported by ANCyT, Secyt-UNC and CIEM}

\dedicatory {To Isabel Dotti and Roberto Miatello on the occasion of their birthday}

\begin{abstract} We prove a Simons-type holonomy theorem for totally skew 1-forms with values in a Lie algebra
of linear isometries. The only transitive case, for this theorem,  is the full orthogonal group. We only use
geometric methods and we do not use any classification (not even that of transitive isometric actions on the
sphere or the list of rank one symmetric spaces). This result was  independently proved, by using an algebraic
approach, by Paul-Andy Nagy. We apply this theorem to prove that the canonical connection of a compact naturally
reductive space is unique, provided the space does not split off, locally, a sphere or a compact Lie group with
a bi-invariant metric. From this it follows easily how to obtain the full isometry group of a naturally
reductive space. This generalizes known classification results of Onishchick, for normal homogeneous spaces with
simple group of isometries, and Shankar, for  homogeneous spaces of positive curvature. This also answers a
question posed by J.\ Wolf and Wang-Ziller. Namely, to explain why  the presentation group of an isotropy
irreducible space, strongly or not, cannot be enlarged (unless for spheres, or for compact simple Lie groups
with a bi-invariant metric).
\end{abstract}

\maketitle

\section{Introduction}

Homogeneous spaces play a central r\^ ole in Riemannian geometry. The most distinguished family is that of
symmetric spaces, defined and classified by E. Cartan \cite{C}. A wider class, that includes the compact
symmetric spaces, are the normal homogeneous ones and, more generally, the  na\-tu\-ra\-lly reductive spaces
\cite{DZ} (see Section 6). Compact isotropy irreducible  spaces, strongly or not \cite{Wo, WZ}, and all known
examples of compact homogeneous Einstein spaces carry naturally reductive metrics.

Symmetric spaces can be defined geometrically, or by means of  a nice presentation, that involves the (connected
component of the) full isometry group. In contrast, the definition of a naturally reductive space $M = G/H$, or
even of a normal homogeneous one, depends on the presentation of the spaces as a quotient of groups. In this case
$G$ needs not to be the full isometry group or its connected component.  For instance, $S^7 = \SO(8) / \SO(7) =
\Spin (7)/G_2$, being  the sphere, with the standard metric,  strongly isotropy irre\-du\-ci\-ble with respect
to both presentations (the same is true for $S^6 = \SO(7)/\SO(6) = G_2/\SU(3)$). The explanation of such a
pathology is that the full isometry group of  this sphere, $\Iso (S^7) = \OO(8)$, does not behave properly with
respect to the canonical connection $\nabla ^c$ associated to the reductive decomposition of $\Spin (7)/G_2$
(since, in general, an isometry maps the canonical connection into another canonical connection).

In general, for a compact naturally reductive homogeneous space $M = G/H$, with a canonical connection $\nabla
^c$, there is a standard way of enlarging the group $G$. Namely, including all the  isometries of $M$ that
belong to $\Aff (M, \nabla^c)$, i.e.\ that are affine with respect to $\nabla ^c$. In fact, if $M$ is simply
connected, any linear isometry $\ell : T_pM \rightarrow T_qM$ with  $\ell(R^c_p) = R^c_q$ and $\ell(T^c_p) =
T^c_q$ extends to an isometry of $M$ (since the canonical connection has $\nabla ^c$-parallel curvature $R^c$
and torsion $T^c$). This extension can also be done, for the connected component and only in the compact case,
by adding to the Lie algebra $\cg$ of $G$, the fields that are invariant under the transvections $\Tr(\nabla
^c)$, which is a transitive normal subgroup of $G$ (such  fields are in a one to one correspondence with the
fixed vectors of the isotropy). This shows that this standard extension is trivial for $S^7 = \Spin(7)/G_2$ (see
Section 7 and \cite{R}).

In this article we shall prove, in a  geometric way, that the sphere is the  only case where this standard
extension does not give the full isometry group.

\begin{teo} \label{fullisometry} Let $M = G/H$ be a naturally reductive Riemannian manifold and let $\nabla^c$
be the associated canonical connection. Assume that $M$ is locally irreducible and that $M \neq S^n$, $M \neq
\bbr P^n$. Then
\begin{enumerate}
\item $\Iso_0(M) = \Aff_0(M, \nabla^c).$
\item If $\Iso (M) \not\subset \Aff (M, \nabla^c)$ then $M$ is isometric to a simple  Lie group, endowed with a
bi-invariant metric (and in this case, the geodesic symmetry maps $\nabla^c$ into $2\nabla - \nabla^c$).
\end{enumerate}
\end{teo}

For a normal homogeneous compact space with $G$ simple, this result was obtained by classification by Onishchik
\cite{On}. Also for homogeneous manifolds of positive curvature the full isometry group was determined by
Shankar \cite{Sh}.

The above theorem follows easily from the following main result which says that the canonical connection is
essentially unique.

\begin{teo} \label{unica} Let $M$ be a compact naturally reductive Riemannian manifold,
which is locally irreducible. Assume furthermore that $M$ is neither (globally) isometric to a sphere, nor to a
real projective space, nor to a compact simple Lie group with a bi-invariant metric.  Then the canonical
connection is unique.
\end{teo}

Taking into account the above result, since the only excluded cases are symmetric spaces,  it make sense to look
for a geometric definition of compact naturally reductive spaces (or at least for the normal homogeneous ones).
Of course, such a definition should coincide with the usual one for  symmetric spaces.

Theorem \ref{unica}  has the following corollary that explains, in a geometric way, that the  presentation of an
isotropy irreducible space gives the connected component of the full isometry group. This question was posed by
J.\ Wolf, for the strongly isotropy irreducible spaces, and by M.\ Wang and W.\ Ziller in general.

\begin{cor}  [\cite{Wo, WZ}] \label{Ziller}  Let $M^n = G/H$ be a compact, simply connected, irreducible homogeneous
Riemannian manifold such that $M$ is not isometric to the sphere $S^n$. Assume that $M$ is isotropy irreducible
with res\-pect to the pair $(G,H)$ (effective action). Assume, furthermore,   that $M$ is not isometric to a
(simple) compact Lie group with a bi-invariant metric. Then $G_0 = \Iso _0(M)$.
\end{cor}

For strongly isotropy irreducible spaces ones needs only to assume that the space is not isometric to a sphere
(see Remark \ref{strongly}).

What it is surprising is that the proof of Theorem \ref{unica}  leads in a natural way to a Simons-type holonomy
theorem that we prove also in a geometric way (involving normal holonomy). In fact, if $\bar \nabla^c$ is
another canonical connection, then the difference tensor at a fixed point $p$, $\Theta = (\bar \nabla ^c -
\nabla ^c)_p$ is a $1$-form with values in the full isotropy algebra (see Section 6). Moreover, since the
Riemannian geodesics, for naturally reductive spaces, are the same as the canonical geodesics, one has that
$\Theta$ is totally skew, i.e.\ $\langle \Theta _X Y, Z\rangle $ is an algebraic $3$-form of $T_pM$ (such a
$1$-form arises naturally as the difference tensor of connections with the same geodesics and are particularly
important in physics, see \cite{AF, A}).

It is then natural to define, just as  Simons \cite{S} did for  holonomy systems, the concept of {\it
skew-torsion holonomy system}: it is a triple $[\bbv, \Theta, G]$, where $\bbv$ is an Euclidean vector space,
$G$ is a connected  Lie subgroup of $\SO(\bbv)$, and $\Theta$ is a totally skew $1$-form on $\bbv$ with values
in the Lie algebra $\cg$ of $G$ (see Section 2).

\begin{teo}[Strong skew-torsion holonomy theorem (see also \cite{N})]\label{SSTHT}
Let $[\bbv, \Theta, G]$, $\Theta \neq 0$, be  an irreducible  skew-torsion holonomy system with $G\neq
\SO(\bbv)$. Then $[\bbv, \Theta, G]$ is symmetric and non-transitive. Moreover,
\begin{enumerate}

\item $(\bbv, [\,\,\,,\,\,])$ is an orthogonal simple Lie algebra, of rank at least $2$,  with respect
to the bracket $[x,y] = \Theta_xy$;
\item $G = \Ad(H)$, where $H$ is the connected Lie group associated to the Lie algebra $(\bbv, [\,\,\,,\,\,])$;
\item $\Theta$ is unique, up to  a scalar multiple.
\end{enumerate}
\end{teo}

This theorem would confirm the hope of J. Simons, in the introduction of \cite{S},  that a variation of his
algebraic setting, for holonomy systems, could be useful for other situations.

In order to apply the above theorem to prove Theorem \ref{unica} one needs a non-trivial result about compact
homogeneous spaces with a nice isotropy group. But this result was also used in the proof of the above theorem
which illustrates the close relationship between both results. Namely,

\begin{teo} \label{Teo3} Let $M = G/H'$ be a compact homogeneous Riemannian ma\-ni\-fold
and let $H$ be the connected component of $H'$. Assume that ($p=[e]$) $T_pM = \bbv_0 \oplus \cdots \oplus
\bbv_k$ (orthogonally) and that $H = H_0 \times \cdots \times H_k$ where
\begin{enumerate}
\item $H_i$ acts trivially on $\bbv_j$ if $i \neq j$ ($i,j = 0, \ldots, k$);
\item If $i \ge 1$, then $H_i$ acts irreducibly on $\bbv_i$. Moreover, $\cc_i(\gh_i) =
\{0\}$, where $\gh_i = \Lie(H_i)$ and $\cc_i(\gh_i) = \{x \in \frak{so}(\bbv_i): [x,\gh_i] = 0\}$.
\end{enumerate}

If $M$ is locally irreducible then $k = 0$ or $k = 1$. Moreover, if $k = 1$ then $\bbv_0 = \{0\}$.
\end{teo}

This result in far from being trivial, due to the fact that  the totally geodesic submanifold $F$, given by the
set of fixed points of the (connected) isotropy $H$, is in general non-trivial. In order to deal with $F$, one
has to use the fact that the $G$-invariant fields are divergence free, which holds in the compact case. If $M$
is not assumed to be compact then Theorem \ref{Teo3} is false; see Remark \ref{contraejemplo}.

The weaker version of Theorem \ref{SSTHT}, where $G$ is assumed to be non-transitive on the sphere, instead of
$G\neq \SO (\bbv)$, follows by using  the same arguments of  the proof of the Simons holonomy theorem \cite{S}
given in \cite{O2}; see Theorem \ref{weakholonomy}.  Such a proof  uses submanifold geometry,  in\-vol\-ving
normal holonomy, as in \cite{O1}. But transitive groups on the sphere, different from the full orthogonal group,
were excluded  by Agricola and Friedrich \cite{AF} (for $\Spin (9)$ after long computations). Actually, the
quaternionic-K\" ahler case was not treated there, but applying the same ideas in \cite{AF}, for the K\" ahler
case, can   be easily solved. This is done  by considering the quaternionic-K\" ahler $4$-form, instead of the
K\" ahler $2$-form (see Remark \ref{quaternionic}).

But such a proof, which involves the classification of  transitive actions on spheres, was far from the spirit
of this article, that  avoids any classification result. We  also succeed in giving a direct and geometric proof
that  transitivity, for skew-torsion holonomy systems, implies  that the group is the full orthogonal group.
This completes the geometric approach to Theorem \ref{SSTHT}.

We would like to remark that many important steps in our arguments relay on submanifold geometry. This
illustrates, as in \cite{O1}, how submanifold geometry can be used to obtain general results in Riemannian
geometry.

This article may be regarded as a step forward in the  geometrization of Lie groups and homogeneous spaces.

We wish to mention that, independently, Paul-Andy Nagy \cite{N}, with an algebraic approach, by means of the
so-called Berger algebras \cite{B}, proved Theorem \ref{SSTHT} (in fact, his preprint appeared earlier than
ours).

\section{Preliminaries and basic facts}

We will extend the definition of Simons to algebraic $1$-forms with values in a Lie algebra, which are totally
skew.

Let $\bbv$ be an Euclidean vector space and let $G$ be  a connected Lie subgroup of $\SO(\bbv)$. Let $\Theta:
\bbv \to \cg = \Lie(G) \subset \mathfrak{so}(\bbv)$ be linear and such that $\langle \Theta_x y, z \rangle$ is
an algebraic $3$-form on $\bbv$. We call such a triple $[\bbv, \Theta, G]$ a {\it skew-torsion holonomy system}.
A skew-torsion holonomy system is said to be:
\begin{itemize}
\item {\it irreducible} if $G$ acts irreducibly on $\bbv$;
\item {\it transitive} if $G$ acts transitively on the unit sphere of $\bbv$;
\item {\it symmetric} if $g(\Theta) = \Theta$, for all $g\in G$,  where
$g(\Theta)_v = g^{-1}\Theta_{g . v} g$.
\end{itemize}

Let $[\bbv, \Theta^\alpha, G]$, $\alpha \in I$, be a family of skew-torsion holonomy systems and let $\cf =
\{g(\Theta^\alpha): g \in G, \alpha \in I\}$. Let $\cg'$  be the linear span of the set $\{\Theta_x: \Theta \in
\cf, x \in \bbv\}$. Then $\cg'$ is an ideal of $\cg$. Let $G'$ be the connected Lie subgroup of $G$ whose Lie
algebra is $\cg'$. Decompose
$$\bbv = \bbv_0 \oplus \cdots \oplus \bbv_k,$$
where $\bbv_0$ is the set of fixed points of $G'$ and $G'$ acts irreducibly on $\bbv_i$ ($i = 1, \ldots, k$).
Let $\cg_i' = \{\Theta_{x_i}: x_i \in \bbv_i, \Theta \in \cf\}$ be the subalgebra of $\cg$, with associated
group $G_i' \subset G$. (cf.\ \cite{S}). Observe that $G_0' = \{0\}$, since $\Theta_{x_0}y = -\Theta_y x_0 = 0$,
for all $x_0 \in \bbv_0$, $y \in \bbv_0^\bot$.

\begin{lema}[see \cite{AF}\label{2.1}, Section 4]
(i) $G' = G_1' \times \cdots \times G_k'$ and $G_i'$ acts irreducibly on $\bbv_i$ and trivially on $\bbv_j$ ($0
\neq i \neq j$). (ii) The decomposition $\bbv = \bbv_0 \oplus \cdots \oplus \bbv_k$ is unique, up to order.
\end{lema}

\begin{proof} Let $\Theta \in \cf$, $x_i \in \bbv_i$, $x_j \in \bbv_j$ and let $x \in \bbv$ ($i \neq j$).
Then $\langle \Theta _xx_i, x_j\rangle = 0 =  \langle \Theta _{x_i}x_j, x\rangle$. Thus $\Theta_{x_i} x_j = 0$
and therefore, if $y = y_1 + \cdots + y_k$, with $y_i \in \bbv_i$, $\Theta_y = \Theta_{y_1} + \cdots +
\Theta_{y_k}$, where $\Theta_{y_i} \in \SO(\bbv_i)$. This implies the lemma.
\end{proof}

\begin{lema}[cf.\ \cite{AF}, Section 4]\label{2.2} Let $\cc_i(\cg_i') = \{x \in \mathfrak{so}(\bbv_i): [x, \cg_i'] = 0\}$.
Then $\cc_i(\cg_i') = \{0\}$ (in particular, $\cg_i'$ is semisimple).

\end{lema}

\begin{proof} Let $B \in \cc _i(\cg_i')$. On has that $\operatorname {ker} (B)$ is $G_i'$-invariant. So, since $G_i'$
acts irreducibly on $\bbv _i$, $\operatorname {ker} (B )= \{0\}$ or $\operatorname {ker} (B) = \bbv _i$. So,
assume that $B\neq 0$ and hence $B$ is invertible. Let $\Theta \in \cf$, $x,y,z \in \bbv _i$. So,
\begin {equation*} \begin {split}
\langle \Theta_x By, z \rangle & = \langle B\Theta_x y, z \rangle =
- \langle B\Theta_y x, z \rangle \\
& = - \langle \Theta_y Bx, z\rangle = \langle \Theta_{ Bx } y,z \rangle.
\end {split}
\end {equation*}

So, $\langle \Theta_{Bx}y,z \rangle = \langle \Theta_x By, z \rangle$. Interchanging $x$ with $z$ we get also
that
$$\langle \Theta_x By, z \rangle = - \langle \Theta_z By, x \rangle = - \langle \Theta_{Bz} y, x \rangle =
\langle \Theta_x y, Bz \rangle.$$ But, since $B$ is skew-symmetric, $\langle \Theta_x y, Bz \rangle = \langle
-B\Theta_x y, z \rangle = - \langle \Theta_x By, z \rangle$. So,
$$\langle \Theta_x By, z \rangle = \langle \Theta_x y, Bz \rangle =
\langle -B\Theta_x y, Bz \rangle = -\langle \Theta_x By, z \rangle.$$ Thus, $\Theta =0$, for all $\Theta \in
\cf$,  a contradiction. Hence $B=0$.
\end{proof}

\begin{prop} \label{2.3} $G = G_0 \times G' = G_0 \times G_1' \times \cdots \times G_k'$, where $G_0$ acts on
$\bbv_0$ and trivially on $\bbv_0^\bot$ ($G_0$ can be arbitrary).
\end{prop}

\begin{proof} We have that $\cg'|_{\bbv_0^\bot}$ is an ideal of $\cg|_{\bbv_0^\bot}$.
Hence, by making use of the previous Lemma, $\cg'|_{\bbv_0^\bot} = \cg|_{\bbv_0^\bot}$. In fact, if $B$ belongs
to the complementary ideal of $\cg'$ in $\cg$, then $B \in \cc(\cg'|_{\bbv_0^\bot})$ and so $B = 0$. From this
it is not hard to prove the proposition, since $\cg'$ is an ideal of $\cg$.
\end{proof}

Let $[\bbv, \Theta, G]$, $\Theta \neq 0$,  be an irreducible skew-torsion holonomy system and let $\nu_v(G . v)$
be the normal space at $v$ to the orbit $G .v$. Namely,
$$\nu_v(G . v) = \{\xi \in \bbv: \langle \xi, \cg . v \rangle = 0\}.$$
We have, as for holonomy systems, see Proposition 3.1 in \cite{O2} , that $\bbw = \nu_v(G . v)$ is a
$\Theta$-invariant subspace, i.e.\ $\Theta_\bbw \bbw \subset \bbw$. In fact, if $\xi \in \bbw$, $x \in \bbv$,
then
$$0 = \langle \Theta_xv, \xi\rangle = -\langle\Theta_\xi v, x\rangle$$
and so $\Theta_\xi . v = 0$. That is $\Theta_\xi \in \cg_v = \Lie(G_v)$ (isotropy algebra). Since $G_v$ leaves
invariant the normal space $\bbw$ of $G . v$ at $v$, one has that $\bbw$ is $\Theta$-invariant.

Then, with the same proof as that given for the Simons holonomy theorem  in \cite{O2} one has:

\begin{teo}  [Weak skew-torsion holonomy theorem] \label{weakholonomy} Let $[\bbv, \Theta, G]$, $\Theta \neq 0$,
be an irreducible non-transitive skew-torsion holonomy system. Then $[\bbv, \Theta, G]$ is symmetric.
\end{teo}

In fact, the proof is slightly more simple since $\Theta$ has less variables than an algebraic Riemannian
curvature tensor.

\begin{prop} \label{symmetricsystem} Let $[\bbv, \Theta, G]$, $\Theta \neq 0$, be an irreducible symmetric skew-torsion
holonomy system. Then:
\begin{enumerate}
\item $G = G'$ (and so,  the linear span of $ \{g(\Theta)_x: g \in G, x \in \bbv\}$ coincides with the Lie algebra
$\cg$ of $G$).
\item $(\bbv, [\,\,\,,\,\,])$ is an (orthogonal) simple Lie algebra
with respect to the bracket $[x,y] = \Theta_xy$;
\item $G = \Ad(H)$, where $H$ is the connected Lie group associated to the Lie algebra
$(\bbv, [\,\,\,,\,\,])$;
\item $\Theta$ is unique, up to a scalar multiple.
\end{enumerate}
\end{prop}

\begin{proof} Part (i) follows from the Proposition \ref{2.3}. If $B \in \cg$ then, since
$[\bbv, \Theta, G]$ is symmetric, $B . \Theta = 0$. So,
$$0 = (B . \Theta)_x y = B \Theta_xy - \Theta_xBy - \Theta_{Bx}y.$$
By making now $B = \Theta_z$,   the Lie identity for the bracket follows. This implies that $\bbv$ is a Lie
algebra. From this it follows (iii). Since $G$ acts irreducibly, it follows that $\bbv$ is simple, which
completes part (ii).

Part (iv) follows from the fact that $\bbv$ is simple. In fact, if $[\bbv, \Theta', G]$ is also a symmetric
skew-torsion holonomy system, then $\Theta'_x$ is a derivation of $(\bbv, [\,\,\,,\,\,])$ and so, $\Theta_x' =
[\ell(x), \,\cdot\;]$, where $\ell: \bbv \to \bbv$ is linear. Since $\Theta'$ and $[\,\,\,,\,\,]$ are
$G$-invariant, we must have that $\ell$ is $G$-invariant, i.e.\ $\ell$ commutes with $G$. Since $G$ acts by
isometries, both the skew-symmetric part $\ell _1$ and the symmetric part $\ell _2$ of $\ell$ must commute with
$G$. By making use of  Lemma \ref{2.2}  we obtain that $\ell _1 = 0$. Moreover, since $G$ acts irreducibly, we
have that $\ell _2 =\lambda \Id$ which proves (iv).
\end{proof}

\begin{nota} \label{rankonegroup} We recall here, for the sake of self-completeness, since we used  this fact,
that there is only one simple connected compact Lie group of rank one (which is the universal cover of $\SO (3)$
and  is isometric to the $3$-sphere). Indeed, let $H$ be a compact simply connected Lie group of rank $1$ and
dimension $n$. Then, as a symmetric pair,  $H = H\times H/\diag(H\times H)$. Since $H$ is a rank-one symmetric
space one has that $\diag (H\times H) \simeq H$ acts transitively on the tangent sphere $S^{n-1}$ at $p = [e]$.
Thus $S^{n-1} = H/S^1$, where $S^1$ is a compact Lie subgroup of dimension $1$ of $H$ (and so, $S^1$ is
homeomorphic to the circle). Recall  that the homotopy groups of $S^1$ are all trivial, except for the first
one. If $n-1 \neq 2$ this yields a contradiction in the exact homotopy sequence induced by $0\longrightarrow S^1
\longrightarrow H\longrightarrow S^{n-1} \longrightarrow 0$. So, $n=3$. In this case the bracket is unique,
since there is a unique $3$-form, up to multiples, in dimension 3. This bracket gives rise to the Lie algebra of
$\SO(3)$.
\end{nota}

\section{The derived 2-form with values in a Lie algebra}

Let $[\bbv, \Theta, G]$, $\Theta \neq 0$, be a skew-torsion holonomy system. We will define a 2-form $\Omega$
with values in $\cg $   such that $\langle \Omega _{x,y}z, w \rangle$ is a $4$-form on $\bbv$.

Let us define
$$\Omega _{x,y} = (\Theta _x . \Theta)_y := [\Theta _x , \Theta _y] - \Theta _{\Theta _x y}$$
It is clear that $\Omega _{x,y} \in \cg$ for all $x, y\in \bbv$. From the definition one obtains  that $\Omega
_{x,y}$ is skew-symmetric in $x$ and $y$. Moreover, for any fixed $x$,  $\langle \Omega _{x,y}z, w \rangle$ is a
$3$-form in the last three variables, since $\Theta$ is totally skew (and then $(B . \Theta)$ is so, for all
$B\in \mathfrak {so}(\bbv)$). Thus $\langle \Omega _{x,y}z, w \rangle$ is an algebraic $4$-form.

\begin{nota} \label{omega} If $v\in \bbv$ then $\Omega _{v,\, \cdot}$ is a (totally skew) $1$-form
with values in the isotropy algebra $\cg _v = \{B \in \cg : B(v)= 0\}$. In fact, $\Omega _{v,\, \cdot} \, v =0
$, since $\Omega$ is totally skew.
\end{nota}

\begin{lema} \label{restrictedform} Let $[\bbv, \Theta, G]$, $\Theta \neq 0$, be  a  skew-torsion holonomy
system and let $\Sigma$ be the  set of fixed points of $H$, where $H$ is a subgroup of $N(G, \OO(\bbv))$ (the
normalizer of $G$ in the full orthogonal group). Assume that the restriction to $\Sigma$ of $\langle \Theta
_{\,\cdot} \, \cdot , \, \cdot \,  \rangle$ is not identically zero. Let $G^\Sigma$ be the connected component
of the subgroup of $G$ that leaves $\Sigma$ invariant. Then:
\begin{enumerate}
\item The cohomogeneity of $G^\Sigma$ on $\Sigma$ is less or equal to the cohomogeneity of $G$ in $\bbv$
(cohomogeneity means the codimension of any principal orbit);
\item There exists a totally skew $1$-form $\Theta ^\Sigma \neq 0$ on $\Sigma$, with values on the Lie algebra
$\cg ^\Sigma$ of $\{g_{|\Sigma} :g\in G^\Sigma\}$ and such that $\langle \Theta^\Sigma _{\,\cdot} \, \cdot , \,
\cdot \,  \rangle$ coincides with the restriction to $\Sigma$ of $\langle \Theta _{\,\cdot} \, \cdot , \, \cdot
\, \rangle$.
\end{enumerate}
\end{lema}

\begin {proof} Part (i) is a special case of Lemma \ref{fixedset}.
Part (ii) follows from the fact that the projection to $\Sigma$ (of the restriction to $\Sigma$) of a Killing
field of $\bbv$, induced by $G$, lies in $\cg ^\Sigma$; see the proof of Lemma \ref{fixedset}. For the sake of
clearness, we adapt the arguments to this particular case. We may assume, by taking the closure, that $H$ is
compact. Let us define
$$\tilde \Theta ^\Sigma = \int_{h \in H} h_*(\Theta).$$
Observe that $\tilde \Theta ^\Sigma$ is a totally skew 1-form with values in $\cg$.

If $ w_1, w_2, w_3 \in \Sigma$, then
\begin{equation*}
\begin{split}
\langle \tilde \Theta^\Sigma_{w_1}w_2, w_3\rangle & = \int_{h \in H} \langle \Theta _{h(w_1)}h(w_2),
h(w_3)\rangle
\\
& = \int_{h \in H} \langle \Theta _{w_1}w_2, w_3\rangle = \langle \Theta _{w_1}w_2, w_3\rangle.
\end{split}
\end{equation*}

Let now $v \in \Sigma^\perp$. Observe that  $\Sigma^\perp$ is $H$-invariant and so, $\int_{h \in H} h(v)= 0$,
since it belongs to  $\Sigma ^\perp$ and it is fixed by $H$. Then

\begin{equation*}
\begin{split}
\langle \tilde \Theta ^\Sigma_{w_1}w_2, v\rangle & = \int_{h \in H} \langle \Theta _{h(w_1)}h(w_2),
h(v)\rangle\\
& = \int_{h \in H} \langle \Theta _{w_1} w_2, v \rangle \\
& = \Big\langle \Theta _{w_1} w_2, \int_{h \in H} v \Big\rangle = 0.
\end{split}
\end{equation*}

Then $\langle \tilde \Theta ^\Sigma_{w_1}w_2, v\rangle =0$. That is, $\tilde \Theta ^\Sigma_{w_1}w_2 \in
\Sigma$. This implies that $\Theta ^\Sigma$, the restriction of $\tilde \Theta ^\Sigma$ to $\Sigma$ has the
desired properties.
\end{proof}

\begin{nota} \label{normal} Let $H$ be a compact Lie group and let $\gh$ be its Lie algebra.
Let $0 \neq v \in \gh$. The normal space to the orbit $H . v : = \Ad (H)v$ is given by
$$\nu_v (H . v) = \cc (v):= \{\xi \in \gh: \ad_v(\xi) = 0\} = \{\xi \in \gh : [v, \xi]=0\}.$$
We have that $v= \Exp (tv) . v = \Ad (\Exp (tv))(v)$ and therefore $\Ad (\Exp (tv))$ leaves $\nu_v (H . v)$
invariant. Moreover, from the above equality,  the set of fixed points  of the one-parameter group of linear
isometries  $\{\Ad (\Exp (tv))\}$ of $\gh$  is just exactly the normal space  $\nu_v (H . v)$.
\end{nota}

\begin{lema}\label{s-representation} Let $[\bbv, \Theta, G]$, $\Theta \neq 0$, be an irreducible
skew-torsion holonomy system. Then $G$ acts on $\bbv$ as an irreducible $s$-representation (i.e. the isotropy
representation of a simple symmetric space).
\end{lema}

\begin{proof} We will define an algebraic (Riemannian) curvature tensor $R\neq 0$ on $\bbv$ with values $R_{v,w}$
in the Lie algebra $\cg$ of $G$. Let
$$R_{v,w} =  [\Theta _v, \Theta _w] - \frac {2}{3} \Omega _{v,w}$$
where $\Omega _{v,w} = (\Theta _v . \Theta)_w$. Then $R_{v,w} \in \cg$ for all $v,w \in \bbv$. Let us verify the
Bianchi identity. Let $\mathcal{B}$ denote the cyclic sum over the first three variables. Since $\Omega _{v,w}z$
is skew-symmetric in $v,w,z$, we have that $\mathcal{B} (\frac {2}{3} \Omega _{v,w}z) = 2\Omega _{v,w}z$. Let us
then compute
\begin{equation*}
\begin{split}
 \mathcal{B} ([\Theta _v, \Theta _w]z) & = [\Theta _v, \Theta _w]z + [\Theta _w, \Theta _z]v +
[\Theta _z, \Theta _v]w \\
& = \Theta _v\Theta _wz - \Theta _{w}\Theta _{v}z + \Theta _{w}\Theta _{z} v - \Theta _{z}\Theta _{w}v
  + \Theta _{z}\Theta _{v}w -\Theta _{v}\Theta _{z}w.
\end{split}
\end{equation*}

Now observe that, in the above sum, the first term is equal to the last one, the second term is equal to the
third one and the remaining two terms are also equal. Thus
\begin{equation*}
\begin{split}
 \mathcal {B} ([\Theta _v, \Theta _w]z) & = 2(\Theta _v\Theta _wz - \Theta _{w}\Theta _{v}z + \Theta _{z}\Theta
 _{v}w)\\
& = 2 [\Theta _v, \Theta _w]z - 2\Theta _{\Theta _{v}w}z  = 2 \Omega _{v,w}z.
\end{split}
\end{equation*}

Then $$\mathcal{B}(R_{v,w}z)= 0.$$

Let us compute $s(R)$, the scalar curvature of $R$. Let $e_1, \ldots , e_n$ be an orthonormal base of $\bbv$.
Since $\langle \Omega _{v,w}z, u\rangle$ is a $4$-form we have
\begin{equation*}
\begin{split}
s(R) & = \sum _{i< j} \langle R_{e_i, e_j}e_j, e_i\rangle = \sum _{i< j} \langle [\Theta_{e_i}, \Theta _ {e_j}]
e_j, e_i\rangle\\
& = \sum _{i< j} \left(\langle \Theta_{e_i} \Theta _ {e_j}e_j, e_i\rangle - \langle \Theta_{e_j} \Theta _
{e_i}e_j, e_i\rangle\right) \\
& = - \sum _{i< j}  \langle \Theta_{e_j} \Theta _ {e_i}e_j, e_i\rangle = \sum _{i< j} \langle \Theta
_{e_i}e_j,\Theta_{e_j} e_i\rangle \\
& = - \sum _{i< j} \langle  \Theta _ {e_i}e_j,\Theta_{e_i} e_j\rangle
\end{split}
\end{equation*}
which is non-zero since $\Theta \neq 0$. Therefore $[\bbv , R, G]$ is an irreducible holonomy system, in the
sense of Simons \cite{S} with $s(R)\neq 0$. Hence $G$ acts as the isotropy representation of a simple symmetric
space by Theorem 5 in \cite{S}. In fact, let $\bar R = \int _{g\in G} g(R)$. Then $s(\bar R) = s(R) \neq 0$ and
therefore $[\bbv , \bar R, G]$ is an irreducible symmetric holonomy system.
\end {proof}

\section{The proof of the strong skew-torsion holonomy theorem}

In this section we state and prove the main result on skew-torsion holonomy systems.

\begin{teo}\label{transitive} Let $[\bbv, \Theta, G]$, $\Theta \neq 0$, be  a transitive  skew-torsion holonomy
system.  Then  $G = \SO(\bbv)$.
\end{teo}

\begin{proof}

Let $\mathcal {F}$ be the family of totally skew $1$-forms on $\bbv$ with values in the Lie algebra $\cg$ of
$G$.  We know, from Section 2, that the linear span of $\{ \tilde \Theta _v: \tilde \Theta \in \mathcal {F}, v\in
\bbv\}$ coincides with $\cg$.

We shall make induction on the dimension $n$ of $\bbv$ ($n \geq 3$, otherwise, $\Theta = 0$). For $n=3$ the
theorem is true, since there is only one such form $\Theta$, up to a scalar multiple, since $\langle \Theta _{\,
\cdot} \, \cdot \, , \, \cdot \,  \rangle$ is a 3-form. In this case  $\{\Theta _v :v\in \bbv \} = \mathfrak
{so}(3)$ and the theorem holds.

Let $n>3$ and assume that the theorem holds for $\dim (\bbv ) < n$.

\smallbreak

{\bf  Case (a).} Assume  that there is   some $\tilde \Theta \in \mathcal {F}$ such that the derived $2$-form
$\tilde {\Omega}$ is not zero, where  $\tilde {\Omega}_{v,w} = (\tilde {\Theta }_v . \tilde {\Theta })_w$.
Choose now some $v\in \bbv$ such that  $\tilde {\Theta} ^v := \tilde {\Omega}_{v,\, \cdot}$ is not zero. We have
that $\tilde {\Theta} ^v$ is a totally skew  $1$-form on $\{v\}^\perp$ with values in the isotropy algebra $\cg
_v$ (see Remark \ref{omega}). So, by the results in Section 2, one has that  the isotropy at $v$, $H=G_v$
satisfies the hypotheses of Theorem \ref{Teo3} with $k\geq 1$, since $\tilde {\Theta} ^v \neq 0$. Observe that,
being the sphere simply connected, one has that the isotropy $G_v$ is connected. Then, since the sphere is
locally irreducible, applying Theorem \ref{Teo3} , one has that  $k=1$ and so,  $H$ acts irreducibly on the
tangent space $T_vS^{n-1}$ of the sphere $S^{n-1} = G . v$.   If the isotropy $H$ is not transitive (on the unit
sphere of the tangent  space $T_vS^{n-1} = \{v \}^\perp $), then $[\{ v \}^\perp, \tilde {\Theta }^v, H]$  is an
irreducible non-transitive skew-torsion holonomy system. Then, by the weak skew-torsion   holonomy Theorem
\ref{weakholonomy}, $[\{ v \}^\perp, \tilde {\Theta }^v, H]$ must be symmetric. Thus we are under the
assumptions of Proposition \ref{grouptype}  to conclude that $S^{n-1} = G . v$ is isometric to a simple Lie
group with a bi-invariant metric which must have rank at least two. This is a  contradiction, since in such a
Lie group, there are totally geodesic and flat submanifolds of dimension at least $2$. Therefore, $H$ is
transitive on the unit sphere of $\{v \}^\perp$. Then $[\{ v \}^\perp, \tilde {\Theta }^v, H]$ is under the
hypotheses of the theorem and $\dim (\{ v \}^\perp) = n-1 < n = \dim (\bbv)$.  So, by induction, $H= \SO(\{ v
\}^\perp)$ which implies that $G= \SO(\bbv)$, since $H=G_v$.

\smallbreak

{\bf Case (b).} Assume that  $\tilde {\Theta} _v . \Theta = 0$ for all $\tilde {\Theta} \in \mathcal {F}$, $v\in
\bbv$ (see Remark \ref{dim4} where it is shown that this case can happen in a non-trivial way).

Observe that, from the assumptions, any $\tilde {\Theta} \in \mathcal {F}$ defines an orthogonal Lie bracket on
$\bbv$. Namely, $[u,v]^{\tilde {\Theta}} := \tilde {\Theta}  _u v$ (and so,  $\ad^{\tilde {\Theta}}_u = \tilde
{\Theta} _u)$.

Let $G^{\tilde {\Theta}}$ be the (connected) Lie subgroup of $G$ associated to the Lie subalgebra $\{\tilde
{\Theta} _v: v\in \bbv\}$ of $\cg$.

By projecting a given  $0\neq \tilde {\Theta} \in \mathcal {F} $ to an irreducible subspace of $G^{\tilde
\Theta}$ (see the decomposition at the beginning of Section 2), we may assume that $\bbv $ decomposes  as $\bbv
= \bbv_0 \oplus \bbv _0^\perp$ into   $\tilde {\Theta}$-invariant subspaces such  that $\tilde {\Theta}$  is
trivial on $\bbv _0$ and irreducible on $\bbv _0^\perp$. Namely, $\tilde {\Theta} _{\bbv _0} =\{0\}$ and
$G^{\tilde {\Theta}}$ acts irreducibly on $\bbv _0^\perp$.

There are three subcases that require different geometric arguments. Name\-ly, $\dim (\bbv_0)\geq 2$, $\dim
(\bbv_0)=1$ and $\dim (\bbv_0)=0$.

For these subcases  we will use that $G$ acts on $\bbv$ as the isotropy representation of a simple symmetric
space (see Lemma \ref{s-representation}). Let us then $0\neq R$ be the unique, up to multiples, algebraic
(Riemannian) curvature tensor on $\bbv$ such that $[\bbv , R, G]$ is a symmetric holonomy system (i.e.\  $g(R) =
R$, for all $g\in G$). In this case $\cg$ coincides with the linear span of $\{R_{u,v}:u,v \in \bbv\}$. We will
show that $R$ has constant curvatures and so $\cg  = \mathfrak {so}(\bbv)$. Since $G$ preserves $R$ and acts
transitively on the unit sphere of $\bbv$, one must only show that there exists $v\neq 0$ such that the Jacobi
operator $J_v = R_{\,\cdot \,  ,\, v}v$,  $J_v : \{v\}^\perp \to \{v\}^\perp$ is a multiple of the identity
transformation.

Before continuing with the proof let us observe the following fact: any normal space at $z$, $\nu _z$ to the
orbit $G^{\tilde {\Theta}} . z$ is $R$-totally geodesic, i.e. invariant under $R$: $R_{\nu _z , \nu _z} \nu _z
\subset \nu _z$. In fact, for any $g\in \{e^{t\tilde {\Theta}_{z}} \}$ and $u, x , y \in \nu _{z}$ one has that,
since $g(R) = R$,
$$ R_{u,x}y = g(R)_{u,x}y = g^{-1} . R_{g . u,g . x} g . y = g^{-1} . R_{u,x}y.$$
So, $R_{u,x}y$ is a fixed vector of   $\{e^{t\tilde {\Theta}_{z}} \}$. Hence, by Remark \ref{normal}, $R_{u,x}y
\in \nu_{z}$ and so, $\nu _{z}$ is left invariant by $R$.

\smallbreak

{\bf (}$\bold{b_1}${\bf ).} Assume that $\dim (\bbv_0 )\geq 2$. Let, for $ v\in \bbv$,   $\nu _v$ denote the
normal space at $v$ of the the orbit $G^{\tilde{\Theta}}.v$. Observe  that   $\bbv _0 \oplus \bbr v \subset \nu
_v$. Observe also that $\nu _v$ is the set of fixed points of the one-parameter group of linear isometries
$\{e^{t\tilde {\Theta}_v} \}$ (see Remark \ref{normal}). There must exist, $v\in \bbv$ and $\bar {\Theta} \in
\mathcal {F}$ such that the associated $3$-form $\langle \bar {\Theta} _{\,\cdot} \, \cdot , \, \cdot \,
\rangle$, when restricted to $\nu _v$, is not identically zero. Otherwise, if $v_0, w_0 \in \bbv _0$ are
linearly independent,
$$\langle \bar {\Theta} _v v_0 , w_0\rangle = 0$$
for all $\bar {\Theta} \in \mathcal {F}$, $v\in \bbv$. Thus, $\bbr v_0 \oplus \bbr w_0$ is perpendicular to any
$G$-orbit, since $\cg$ is linearly spanned by $\{\bar {\Theta}_v:\bar {\Theta } \in \mathcal {F}, v\in \bbv \}$.
This is a contradiction, since $G$ is transitive on the sphere.

Let us then choose  $v \in \bbv$ and $ \bar {\Theta}$  such that  its associated $3$-form is not identically
zero when restricted to $\nu_v$. By perturbating $v$ slightly, we may assume that $v\notin \bbv_0$. Moreover, if
$v'$ is the orthogonal projection of $v$ to $\bbv_0 ^\perp$, one has that $\nu _{v'} = \nu _v$. So, we may
assume that $0\neq v \in \bbv _0^\perp$.

We can now apply Lemma \ref{restrictedform} to conclude that the cohomogeneity of $G^{\nu_v}$ on $\nu _v$ is one
(i.e.\ it is transitive on the sphere). Moreover, there exists a totally skew $1$-form $ \bar {\Theta }^v \neq
0$ with values in the Lie algebra  $\cg^v$ of $G^{\nu_v}$ (keeping the notation of Lemma \ref{restrictedform}).
Since $\dim (\nu _v) < n = \dim (\bbv)$ we have, by induction, that $\{g_{|\nu_v}: g\in G^{\nu_v}\} = \SO
(\nu_v)$. Let now $v_0$ and $w_0$ be perpendicular vectors of unit length that belong both to $\bbv _0$ and let
$$\lambda = \langle J_{v_0}(w_0), w_0\rangle.$$

Let now $z\in \bbv$ be an arbitrary vector which is perpendicular to $v_0$. Let $v' = g(v)$, $g\in G^{\tilde
{\Theta }}$, be such that $z\in \nu _{v'}$ (e.g.\ by choosing $v' \in G^{\tilde {\Theta}} . v$ such that the
high function $x \mapsto \langle z, x\rangle$ on this orbit attains its maximum value). Since $\nu _{v'} = g(\nu
_v)$ and $G^{\nu _{v'}} = g G^{\nu _{v}}g^{-1}$, we have that $G^{\nu _{v'}} = \SO(\nu _{v'})$. Recall that $R$
leaves invariant   the subspace $\nu _{v'}$ (that contains the vectors $v_0,w_0$). Since $G^{\nu _{v'}} =
\SO(\nu _{v'})$ preserves $R$,  we conclude that the restriction of $R$ to $\nu _{v'}$ has constant sectional
curvatures $\lambda$. So, in particular $J_{v_0}(z) = \lambda z$. Thus $J_{v_0} : \{v_0\}^\perp \to
\{v_0\}^\perp$ is a multiple of the identity. This concludes the proof of this subcase.

\smallbreak

{\bf (}$\bold{b_2}${\bf ).} Assume that $\dim (\bbv_0 ) = 1$. Let $v\in \bbv_0$ be of unit length. Observe that
$G^{\tilde {\Theta}}$ preserves $R$, fixes $v$ and acts irreducibly on $\{v\}^\perp$. Then $G^{\tilde {\Theta}}$
commutes with $J_v$ and therefore $J_v$ is a multiple of the identity.

\smallbreak

{\bf (}$\bold{b_3}${\bf ).} Assume that $\dim (\bbv_0 ) = 0$, i.e.\ $G^{\tilde {\Theta}}$ acts irreducibly on
$\bbv$. In this case the principal orbits of $G^{\tilde {\Theta}}$ are irreducible and full isoparametric
submanifolds of $\bbv$. In fact,  $G^{\tilde {\Theta}}$ acts as the $\Ad$-representation of $(\bbv ,\,  [\ , \
]^{\tilde {\Theta}})$ (see \cite{PT}, \cite{BCO}). The cohomogeneity of  $G^{\tilde {\Theta}}$ in $\bbv$ is at
least $2$. Otherwise, the Lie group associated to $(\bbv ,\,  [\ , \ ]^{\tilde {\Theta}})$ would be of rank one
and therefore, by Remark \ref{rankonegroup}, $\dim (\bbv) = 3$, a contradiction, since we assume $n>3$.

Let $M = G^{\tilde {\Theta}} . v$  be a principal orbit. Let $\xi$ belong to the normal space $\nu_v$ of $M$ at
$v$  such that the shape operator $A_\xi$ has all of  its eigenvalues $\lambda _1, \ldots , \lambda _g$
different from zero and $g\geq 2$. Such a $\xi$ can be chosen by perturbating slightly the position vector,
since the the codimension of $M$ is at least $2$ and so $M$, since it is full,  is not umbilical. Let $E_1,
\ldots, E_g$ be the eigenspaces of $A_\xi$ associated to $\lambda _1, \ldots , \lambda _g$ respectively. Let us
write
$$\bbv = \nu _v \oplus E_1 \oplus \cdots  \oplus E_g.$$
Let  $v_i  = v + \lambda _i^{-1}\xi$ ,  $i=1, \ldots , g$. One has that the normal space at $v_i$ to the orbit
$M_i = G^{\tilde {\Theta}} . v_i$ is
$$\nu_{v_i}= \{ \xi \in \bbv : \tilde {\Theta }_{\xi}v_i  = 0\} =  \nu _v \oplus E_i$$
which is $\tilde {\Theta }$-invariant, since $\{ \xi \in \bbv : \tilde {\Theta }_{\xi}v_i  = 0\}$ is a
subalgebra of $(\bbv ,\,  [\ , \ ]^{\tilde {\Theta}})$. Moreover, the restriction of $\tilde {\Theta }$ to $\nu
_{v_i}$ is not zero, since $\nu _{v_i}$ is not abelian. In fact,  it contains properly the maximal abelian
subalgebra $\nu _v$. Recall that, by Remark \ref{normal},  $\nu _{v_i}$ is the set of fixed points of
$\{e^{t\tilde {\Theta}_{v_i}} \}$ (and so, as we previously remarked, $\nu _{v_i}$ is left invariant by $R$). By
Lemma  \ref{restrictedform}, $G^{\nu _{v_i}}$ is transitive on the unit sphere of $\nu _{v_i}$. Hence, by
induction, since $\dim (\nu _{v_i}) < \dim (\bbv)$,  we get that $G^{\nu _{v_i}} = \SO(\nu _{v_i})$.

Since the restriction $R^i$ of $R$ to $\nu_{v_i}$ is fixed by $G^{\nu _{v_i}} = \SO(\nu _{v_i})$, one has that
$R^i$  has constant curvatures, let us say $\mu$. Let $\bbw _i = \{v\}^\perp \cap \nu _i$ then
$$J_v{_{|\bbw _i}} = \mu \Id _{\bbw _i}$$

We have  dropped the subscript $i$ for $\mu$,  since it does not depend on $i$. In fact, let  $w \in \nu _v$ be
of unit length and  perpendicular to $v$.

We have that  $w$ as well as $v$,  belong both to any $\bbw _i$ and that $\mu = \langle R_{w,v}v, w \rangle$.
This shows that $\mu$ is independent of  $i =1, \ldots , g$. Since, $ \{v\}^\perp = \bigcup _{i} \bbw _i$ we
conclude that $J_v$  coincides on $ \{v\}^\perp$ with $\mu \Id_ {\{v\}^\perp}$.

\smallbreak

This concludes the proof of the theorem.
\end{proof}

\begin{proof}[Proof of Theorem \ref{SSTHT}]  By putting together the  weak skew-torsion holonomy Theorem
\ref{weakholonomy}, Proposition \ref{symmetricsystem} and  Theorem  \ref{transitive}, Theorem \ref{SSTHT} follows.
\end{proof}

\begin{nota}\label{dim4}
Let us consider in $\bbr ^4 \simeq \bbr \oplus \mathfrak {so }(3)$ the bracket given by the product of  the
(trivial) bracket of $\bbr$ and the bracket  of $\mathfrak {so}(3)$. Any bracket in $\bbr ^4$ defines in a
natural way a $3$-form. Since the space of $3$-forms in $\bbr ^4$ is canonically isometric to $\bbr ^4$, then
the group $\SO(4)$ acts transitively on the family of $3$-forms of unit length. This implies  that any $3$-form
defines a bracket in $\bbr ^4$, which is orthogonally equivalent to the given one (up to a scalar multiple). Let
now  $ \Theta \neq 0$ be any   totally skew $1$-form with values in $\mathfrak {so}(4)$. Then $\Theta$ satisfies
the equation
$$\Theta _x . \Theta  = 0 $$
for any $x\in \bbr ^4$. But $[\bbr ^4, \Theta , \SO (4)]$ is never a symmetric skew-torsion holonomy system.
\end{nota}

\begin{nota} \label{quaternionic}

By using  the classification of the  transitive actions on spheres one can prove Theorem \ref{transitive}. In
fact, all the cases have been excluded by Agricola and Friedrich \cite{AF} except for the quaternionic-Kh\"aler
case $G= \Sp(1)\times \Sp(n)$ with $n>1$. But just copying their argument, used for the K\" ahler case, one can
also exclude this group. Instead of using the K\" ahler form, one has to use the quaternionic-K\" ahler
$4$-form.

Let
$$ \Omega = e_1\wedge e_2 \wedge e_3 \wedge e_4 + \dots + e_{4n-3}\wedge e_{4n-2} \wedge e_{4n-1} \wedge
e_{4n}$$ be the quaternionic-K\" ahler  4-form which is left invariant by $\Sp(1)\times \Sp(n)$ (nicely written
in an appropriate basis). Let $\Theta$ be a totally skew $1$-form with values in $\mathfrak {sp}(1)\times
\mathfrak {sp}(n)$ and let  $\Theta  _{i,j}^k = \langle \Theta _{e_i}e_j, e_k\rangle$. Since $\Omega $ is left
invariant by $\Sp(1)\times \Sp(n)$, we must have that $\Theta _v . \Omega = 0 $ for all $v \in \bbr ^{4n}$.

Let $i,j,k$ be different indexes  $\{i,j,k\} \subset \{1,2,3,4\}$. We assume that $i = 1$, $j=2$, $k=3$ (the
other cases are similar). Let $m > 4$ and $r$ be arbitrary.
\begin{equation*}
\begin{split}
0 & = (\Theta _{e_r}. \Omega)(e_1, e_2 , e_3 , e_m) \\
& = \Omega (\Theta _{e_r} e_1 , e_2, e_3, e_m) + \Omega (e_1, \Theta _{e_r} e_2 , e_3, e_4)\\
& \qquad + \Omega ( e_1 , e_2, \Theta _{e_r}e_3, e_m) + \Omega ( e_1 , e_2, e_3, \Theta _{e_r} e_m) \\
&= \Theta _{r,m}^4 \Omega(e_1, e_2 , e_3 , e_4) = \Theta _{r,m}^4 = - \Theta _{r,4}^m.
\end{split}
\end{equation*}

So, $\Theta _{r,4}^m = 0$. Similarly, $\Theta _{r,3}^m = \Theta _{r,2}^m = \Theta _{r,1}^m = 0 $. The same is
true if we replace $\Theta $ by $g(\Theta)$, for all $g \in \Sp(1)\times \Sp(n)$.

Then, for all $g \in \Sp(1)\times \Sp(n)$ and for all $m > 1$,
$$g(\Theta) _{r,1}^m =  g(\Theta) _{r,2}^m = g(\Theta) _{r,3}^m  = g(\Theta) _{r,4}^m  = 0.$$
Observe that Lemma \ref{2.1} and Proposition \ref{2.1} imply that $\mathfrak {sp}(1)\times \mathfrak {sp}(n)$ is
linearly spanned by $\{g(\Theta)_{v}: v \in \bbr ^{4n}\}$. Then, the above equality implies that $\Sp(1)\times
\Sp(n)$ leaves invariant the subspace of $\bbr ^{4n}$ spanned by $\{e_1,e_2, e_3, e_4\}$. This is a
contradiction since $\Sp(1)\times \Sp(n)$ acts irreducibly on $\bbr ^ {4n}$ ($n>1)$. Therefore, $\Sp(1)\times
\Sp(n)$ is also excluded.
\end{nota}

\section{Compact homogeneous manifolds with nice isotropy}

In this section we prove a general result, Theorem \ref{Teo3}, for compact homogeneous Riemannian manifolds with
nice isotropy groups.

Let us first state some  results that we need.

\begin{lema} \label{normaldistribution} Let $M = G/G_p$ be a homogeneous Riemannian manifold, let $H$
be a normal subgroup of $G_p$ and let $\bbw$ be the subspace of $T_pM$ defined by
$$\bbw = \{v \in T_pM: dh|_p(v) = v, \text{ for all  } h \in H\}.$$
Then $\bbw$ is $G_p$-invariant. Moreover, if $\cd$ is the $G$-invariant distribution on $M$ defined by $\cd_p =
\bbw$, then $\cd$ is integrable with totally geodesics integral manifolds (or, equivalently, $\cd$ is
autoparallel).
\end{lema}

\begin{proof} Let us construct explicitly the integral manifold $S(q)$ by $q$ of $\cd$. Let $q = g . p$ and
let
$$S(q) = \{x \in M: h . x = x, \text{ for all }h \in gHg^{-1}\}$$
be the fixed points of $gHg^{-1}$. It is well known that $S(q)$ is a totally geodesic submanifold of $M$ (we may
assume that $H$ is compact, otherwise the closure of $H$ is also a normal subgroup of $G_p$ with the same fixed
vectors and fixed set). One has that $T_qS(q) = \cd_q$. Let now $r \in S(q)$. Since $S(q)$ is a homogeneous
submanifold of $M$ (see Lemma \ref{fixedset}) there exist $g' \in G$ with $g' . q = r$ and such that $g'(S(q)) =
S(q)$. Then $T_rS(q) = dg'(T_qS(q)) = dg'(\cd_q) = \cd_r$, since $\cd$ is $G$-invariant. Then $S(q)$ is an
integral manifold of $\cd$.
\end{proof}

\begin{lema} \label{fixedset} Let $M$ be a  Riemannian manifold, $G$ a closed connected subgroup of
the isometries $\Iso (M)$ of $M$  and let $H \subset N(G,\Iso(M))$ (the normalizer of $G$ in  the full isometry
group). Let
$$\Sigma = \{x \in M: h(x) = x \text{ for all } h \in H\}$$ be the set of fixed points of $H$,
which we assume to be non-empty (observe that $\Sigma$ is a closed and totally geodesic submanifold of $M$). Let
$G^\Sigma$ be the connected component of the subgroup of $G$ that leaves $\Sigma$ invariant. Then the
cohomogeneity of $G^\Sigma$ on $\Sigma$ is less or equals to the cohomogeneity of $G$ on $M$. In particular, if
$G$ is transitive on $M$, then $G^\Sigma$ is transitive on $\Sigma$.
\end {lema}

\begin{proof} We may assume that $H$ is a closed Lie subgroup. The group $H$ turns out to be compact since
any element of $\Sigma$ is a fixed point of $H$. Let us endow $H$ with an $H$-invariant
volume element such that $ \vol(H) = 1$. Let $X \in \ck _G (M) \simeq \cg$ (Killing fields of $M$ induced by
$G$) and define $\bar X \in \ck _G(M)$ by
$$\bar X = \int_{h \in H} h_*(X) \ \  \in \, \ck _G (M).$$

We have that $\bar X_r$ is the projection to $T_r\Sigma$ of $X_r$, for all $r \in \Sigma$. In fact,
\begin{equation*}\begin{split}
\bar X_r & = \int_{h \in H} dh|_{h^{-1}(r)}(X_{h^{-1}(r)}) \\
& = \int_{h \in H} dh|_r(X_r)\\
& = \int_{h \in H} dh|_r(v) + \int_{h \in H} dh|_r(w),
\end{split}\end{equation*}
where $X_r = v + w$, $v \in T_r\Sigma$, $w \in (T_r\Sigma)^\bot$. One has that
$$\int_{h \in H} dh|_r(v) = \int_{h \in H} v = v.$$
On the other hand, the vector $z = \int_{h \in H} dh|_r(w)$ is perpendicular to $T_r\Sigma$ and it is fixed by
$H$. This implies that $z = 0$. So, $\bar X|_\Sigma$ is always tangent to $\Sigma$. Moreover, $\bar X|_\Sigma$
coincides with the projection to $T\Sigma$ of $ X|_\Sigma$.

Let $\ck _G (\Sigma)$ be the Killing fields of $M$, induced by $G$, and such that, when restricted to $\Sigma$,
are always tangent to $\Sigma$. Then, $\ck _G (\Sigma)$, coincide with the projection to $\Sigma$, of the
restriction to $\Sigma$, of the elements of $\ck _G (M)$. It is now clear that a vector in $T_r\Sigma$ which is
perpendicular to the  orbit $G^\Sigma . r$ must be perpendicular to $G . r$. This implies the lemma.
\end{proof}

The following lemma is crucial for our purposes.

\begin{lema} \label{3.3} Let $M = G/H$ be a compact homogeneous Riemannian ma\-ni\-fold, where $H$ is connected. Let
$\bbv_0 \subset T_pM$ be the set of fixed vectors at $p = [e]$ of $H$. Assume that $H$ acts irreducibly on
$\bbv_0^\bot$ and that $\cc(\gh) = 0$, where  $\gh = \Lie(H)$ and $\cc(\gh) = \{B \in \frak{so}(\bbv_0^\bot):
[B, \gh] = 0\}$.
 Let $\cd$ be the $G$-invariant distribution on $M$ with $\cd_p = \bbv_0$. Then $\cd$ is a parallel
distribution on $M$ (and hence, if $\cd$ in non-trivial, $M$ splits locally).
\end{lema}

\begin{proof} Let $v \in \bbv_0$ and let $\tilde v$ be the $G$-invariant vector field on $M$ with $\tilde v_p = v$.
Since $M$ is compact, $\tilde v$ divergence free on $M$ (since the flow associated to $\tilde v$ commutes with
$G$ and so it must preserve volumes; see Remark \ref{volume}). Moreover, if $S(q)$ is an integral manifold of
$\cd$ though $q$ (which is totally geodesic in $M$) then $\tilde v|_{S(q)}$ is also divergence free in $S(q)$.
In fact, let $\tilde G(q) = \{g \in G: g(S(q)) = S(q)\}$. Then $\tilde G(q) \subset G$ is a closed Lie subgroup
which acts transitively on the compact manifold $S(q)$ (see Lemma \ref{fixedset}). Then $\tilde v|_{S(q)}$ is
$\tilde G(q)$-invariant. Hence $\tilde v_{S(q)}$ is divergence free in $S(q)$.

 Since $\cd$ is $G$-invariant, it
suffices to show that $\cd$ is parallel at $p$. If $h \in H$ is arbitrary,
$$dh|_p(\nabla_u \tilde v) = \nabla_{dh|_p(u)} h_*(\tilde v) = \nabla_{dh|_p(u)} \tilde v,$$
for all $u \in T_pM$. Then, if $T:T_pM \to T_pM$, $T(u) = \nabla_u\tilde v$, $T$ commutes with $H$ (via the
isotropy representation). Then $H$ commutes with both, the symmetric part, let us say $A$, and the
skew-symmetric part, let us say $B$, of $T$. In particular, both $A$ and $B$ leave $\bbv_0$ and $\bbv_0^\bot$
invariant. Since $B|_{\bbv_0^\bot} \in \cc(\gh)$ then, from the assumption, $B|_{\bbv_0^\bot} = 0$.

Being  $\tilde v$  divergence free on $M$, one has that $\trace (A) = 0$. Since $\tilde v|_{S(p)}$ is also
divergence free, $\trace(A|_{\bbv_0}) = 0$ and therefore $\trace(A|_{\bbv_0^\bot}) =  0$. Since $H$ acts
irreducibly on $\bbv_0^\bot$, $A|_{\bbv_0^\bot} = \lambda \Id$. Thus,   $A|_{\bbv_0^\bot} = 0$.

Then $(\nabla_{T_pM}\tilde v)_p \subset \bbv_0 =  \cd_p$. Then $\cd$ is parallel at $p$. Hence, since $\cd$ is
$G$-invariant, $\cd$ is parallel on $M$.
\end{proof}

\medbreak

\begin{proof}[Proof of Theorem \ref{Teo3}] We may assume that $G$ is compact. In fact, let $(\bar G_p)_0$ be the
connected component of the isotropy group of the closure $\bar G$ of $G$ in $\Iso (M)$, which must belong to
the normalizer of $H$. Then $(\bar G_p)_{0|\bbv _i}= H_{|\bbv _i}$, $i=1, \ldots , k$. Otherwise, any $0\neq X$
in a complementary ideal of $\gg _i$ (inside the Lie algebra  of $(\bar G_p)_{0|\bbv _i}$) would belong to
$\cc_i(\gh_i)$. This implies that $(\bar G_p)_0 = \tilde H_0 \times H_1 \times \cdots \times H_k$, where $\tilde
H_0$ acts only on $\bbv _0$. Then $(\bar G_p)_0$ is also in the assumptions of the theorem.
  We also may assume, by passing to a finite cover,
that $H' = H$.

The key fact is to prove that the distribution given by the fixed vectors of the isotropies is parallel along
the $G$-invariant distributions defined by $\bbv_0^\bot$.

Assume that $k \ge 1$ and let
$$H^1 = H_0 \times H_2 \times \cdots \times H_k$$
which is, as well as $H_1$, a normal subgroup of $H$. Let $\bbv^1 = \bbw_0 \oplus \bbv_1$ be the set of fixed
vectors at $p$ of $H^1$, where $\bbw_0 \subset \bbv_0$ is the set of fixed vectors at $p$ of $H_0$. The subspace
$\bbv^1$ is $H$-invariant and so it extends to a $G$-invariant distribution $\cd^1$ on $M$, which is also
autoparallel (see Lemma \ref{normaldistribution}). Observe that the integral manifold $S^1(p)$ by $p$ is the
(connected component containing $p$ of) the set of fixed points of $H^1$ in $M$. We have that $S^1(p)$ is
homogeneous under the group (see the previous lemma)
$$G^1(p) =  \{g \in G: g(S^1(p)) = S^1(p) \} = \{g \in G: g(p) \in S^1(p) \}.$$
Observe that the isotropy $(G^1(p))_p$ coincides with $H$. But $H^1$ acts trivially on $S^1(p)$. Then the
effectively made isotropy of $S^1(p)$ is $H_1$. We have that $\bbw_0 \subset T_pS^1(p)$ is the set of fixed
vectors at $p$ of $H_1$ in $S^1(p)$. Then,  by the previous lemma,  the $G$-invariant distribution $\cd_0$, with
$\cd_0(p) = \bbw_0$ is parallel along $S^1(p)$. Thus, the $G$-invariant distribution $\cd _1$ on $M$, with $\cd
_1(p) = \bbv_1$ is also parallel along $S^1(p)$. By the $G$-invariance of $\cd_1$ and $\cd^1$, we also have that
$\cd_ 1$ is parallel along $S^1(q)$, for all $q \in M$. This implies that $\cd_1$ is an autoparallel
distribution on $M$ (since, $\cd ^1$ is autoparallel and $\cd _1 \subset \cd ^1$). But $\cd_1^\bot$ is also an
autoparallel distribution on $M$. In fact, $\cd_1^\bot$ is the $G$-invariant distribution with $(\cd_1^\bot)_p =
\bbv_0\oplus \bbv_2 \oplus \cdots \oplus \bbv_k$, which is the set of fixed vectors at $p$ of $H_1$. But two
orthogonally complementary autoparallel distributions must be parallel (see e.g.\ \cite[p.\ 31]{BCO}). Then $M$
splits locally, unless $\cd_1^\bot = \{0\}$.
\end{proof}

\begin{nota} \label{contraejemplo}
Theorem \ref{Teo3} does not hold if $M$ is not assumed to be compact. In fact, let $H^n$, $n \ge 4$, be the real
hyperbolic space and let $\mathcal{F}$ be a foliation of $H^n$ given by parallel horospheres, centered at same
point $q_\infty$ at infinity. Let $G$ be the (connected component of) the subgroup of $\Iso_0(H^n)=
\SO_0(n+1,1)$ that leaves $\mathcal{F}$ invariant. Then $G$ acts transitively on the hyperbolic space, since it
contains the solvable group that fixes the point $q_\infty$. Let $p \in H^n$ and let $v\in T_pM$ be
perpendicular to the horosphere of $\mathcal {F}$ by $p$. Then the isotropy $G_p$, via the isotropy
representation, fixes $v$. Moreover, $G_p$ restricted to $\{v\}^\perp$ coincides with $\SO (\{v\}^\perp) \simeq
\SO(n-1)$. If Theorem \ref{Teo3} holds, then $H^n$ would be reducible (in this case it would split off a line),
which is a contradiction.
\end{nota}

\begin{nota}\label{volume}
Let $M$ be a Riemannian manifold and let $G$ be a Lie group acting transitively on $M$ by isometries and such
that $G$ admits a bi-invariant metric $(\,\, ,\, )$. Let $\varphi$ be a difeomorphism of $M$  that commutes with
$G$. Then, as it is well known,  $\varphi$ is volume preserving. In fact, let $B_\varepsilon $ be an
$\varepsilon$-ball,with respect to $(\,\, ,\, )$, around the identity $e\in G$. If  $p\in M$, then
$\{B_\varepsilon .p\}$ is a basis of neighborhoods of $p$ in $M$, for arbitrary $\varepsilon >0$. Let $g\in G$
be such that $g.p = \varphi (p)$. Then
$$\varphi(B_\varepsilon .p) = B_\varepsilon .  \varphi(p) =   B_\varepsilon . (g.p)
= g.(g^{-1}B_\varepsilon g).p = g. B_\varepsilon .p $$ and so,
$$\text {vol} (\varphi(B_\varepsilon .p)) = \text {vol} (g. B_\varepsilon .p)
=  \text {vol} ( B_\varepsilon .p)$$ This implies that $\varphi$ is volume preserving.
\end{nota}

\section{Applications to naturally reductive spaces}

Let $M = G/H$ be a homogeneous compact Riemannian manifold with a $G$-invariant metric $\langle \;\,,
\;\rangle$. The space $M$ is said to be {\it naturally reductive} if there exists a reductive decomposition
$$\cg = \gh \oplus \gm,$$
where $\cg = \Lie(G)$, $\gh = \Lie(H)$, $\Ad(H)\gm \subset \gm$, such that the geodesics by $p = [e]$ are given
by
$$\gamma_{X . p} = \Exp(tX) . p$$
for al $X \in \gm$. In other words, the Riemannian geodesics coincide with the $\nabla^c$-geodesics, where
$\nabla^c$ is the canonical connection, which is a metric connection,  of $M$ associated to the reductive
decomposition. This is in fact equivalent to the property that $[X, \,\cdot\;]_\gm : \gm \to \gm$ is
skew-symmetric, for all $X \in \gm$ ($\gm \simeq T_pM$).

The Levi-Civita connection is given by
$$\nabla_v \tilde w = \mbox{$\frac{1}{2}$} [\tilde v, \tilde w]_p,$$
and
$$\nabla_v^c \tilde w = [\tilde v, \tilde w]_p,$$
where, for $u \in T_pM$, $\tilde u$ is the Killing field on $M$ induced by the unique $X \in \gm$ such that $X .
p = u$ (i.e.\,  $\tilde u(q) = X . q$).

The difference tensor between both connections is given by
$$D_vw = \nabla_v \tilde w - \nabla_v^c \tilde w = -\mbox{$\frac{1}{2}$} [\tilde v, \tilde w]_p
= -\nabla_v \tilde w.$$ The tensor $D$ is totally skew, i.e.\ $\langle D_vw, z \rangle$ is a $3$-form.

Assume that the Riemannian metric of $M$ is also naturally reductive with respect to another decomposition. That
is, $M = G'/H'$ and there is a reductive decomposition
$$\cg ' = \gh ' \oplus \gm ',$$
such that the geodesics of  $M$  by $p$ are also given by
$$\gamma_{X . p} = \Exp(tX) . p$$
for $X\in \gm '$, where $\cg '  = \Lie(G')$, $\gh '  = \Lie(H')$, $\Ad(H')\gm '  \subset \gm '$,

Let $\bar \nabla ^c$  be the canonical connection associated to this new reductive decomposition. Similarly, if
$\bar D  = \nabla - \bar \nabla ^c$,
$$\bar D_vw = \nabla_v \bar w - \bar {\nabla}_v^c \bar w = -\mbox{$\frac{1}{2}$} [\bar v, \bar w]_p
= -\nabla_v \bar w.$$ where, for $u \in T_pM$, $\bar u$ is the Killing field on $M$ induced by the unique $X \in
\gm '$ such that $X . p = u$ (i.e.\, $\bar u(q) = X . q$).

The  tensor $\bar D$, as well as $D$ is totally skew.

One has that
$$D_vw - \bar D_vw = -\nabla_v( \tilde w - \bar w) = - \nabla_v Z,$$
where $ Z = \tilde w - \bar w$ vanishes at $p$. So, $(\nabla Z)_p \in \tilde \gh := \Lie(\Iso(M)_p)$ (via the
isotropy representation). In fact,
$$e^{t(\nabla Z)_p} = d (\varphi_t^Z)|_p,$$
where $\varphi_t^Z$ is the flow associated to $Z$. Observe that $D - \bar D = \bar \nabla^c - \nabla^c$.

Then, if $\Theta : = D - \bar D$, $\Theta w = -(\nabla Z)_p \in \tilde \gh$  or, equivalently,
$$\Theta _w = (\nabla Z)_p \, \in \tilde \gh$$
since $\Theta$ is totally skew. That is, $\Theta _w$ belongs to the full isotropy algebra, for all  $w \in
T_pM$.

Let
$$\bar \gh = \text {linear span of $\{g(\Theta )_w: g \in \Iso(M)_p,\, w \in T_pM\}$}.$$
One has, as in Section 2, that $\bar \gh$ is an ideal of $\tilde \gh = \Lie(\Iso(M)_p) \subset \frak{so}(T_pM)$.

Let $\bar H$ be the connected Lie subgroup of $\SO(T_pM)$ with $\Lie(\bar H) = \bar \gh$. Then, by what has been
done for skew-torsion holonomy systems (see Section 2),
$$T_p M = \bbv_0 \oplus \bbv_1 \oplus \cdots \oplus \bbv_k \quad \text{(orthogonally)}$$
and
$$(\Iso_0(M)_p)_0 = H_0 \times H_1 \times \cdots \times H_k,$$
where $H_0$ acts only on $\bbv_0$ and $H_i$ acts irreducibly on $\bbv_i$ and trivially on $\bbv_j$ if $i \neq
j$, $i \ge 1$. Moreover, such groups are in the assumptions of Theorem \ref{Teo3}. Then, by this theorem, if $M$
is locally irreducible, either
$$\text{$(\Iso_0(M)_p)_0 = H_0$ and $T_pM = \bbv_0$},$$
or
$$\text{$(\Iso_0(M)_p)_0 = H_1$ and $T_pM = \bbv_1$}.$$

If $(\Iso_0(M)_p)_0 = H_0$ then all $ g(\Theta) = 0$, and in particular $\Theta = 0$ and so $\nabla ^c = \bar
\nabla ^c$.

Let us analyze the remaining case, where $(\Iso_0(M))_0 = H_1$.

From  Theorem \ref{SSTHT}, there are only two cases:

{\bf (a)} $H_1 = \SO(T_pM)$. In this case $M$ has (positive) constant curvatures and it is a global symmetric
space (see Remark \ref{projective}). Then $M = S^n$ or $M = \bbr P^n$.

{\bf (b)} $H_1$ acts on $T_pM$ as the $\Ad$-representation of a compact simple Lie group. Then $M$ is isometric
to a compact simple group with a bi-invariant metric (a classification free and geometric proof of this fact is
given in Proposition \ref{grouptype})

\begin{proof}[Proof of Theorem \ref{unica}] It follows from the above discussion.
\end{proof}

\begin {nota} \label{pruebaaffinecompact} We keep the above notation. Assume, as in case (b),
that $M$ is isometric to a compact simple Lie group with a bi-invariant metric. In this case the family of
canonical connections is the affine line
$$L = \{s\nabla + (1-s)\nabla^c:s \in \bbr\},$$
since the difference tensor between any two canonical connections must be, up to a scalar multiple, unique (see
Proposition \ref{symmetricsystem}). We assume $\nabla^c \neq \nabla$. In this case $\Iso_0(M)$ fixes point-wise
the whole line (since $\Iso(M)$ induces an isometry on $L$ with the fixed point $\nabla$). In particular $\Iso
_0(M)$ preserves $\nabla ^c$. So, $\Iso _0 (M) \subset \Aff_0(M, \nabla^c)$.  Observe that in this case the
geodesic symmetry moves the totally skew tensor $D$ into $-D$. So, it moves $\nabla ^c = \nabla - D$ into
$\nabla + D = 2\nabla + (-\nabla + D) = 2\nabla - \nabla^c$. Then the geodesic symmetry does not belong to $\Aff
(M, \nabla ^c)$.
\end {nota}

\begin{proof}[Proof of Theorem \ref{fullisometry}] Assume that  $M$ is not also isometric to a compact (simple)
Lie group with with a bi-invariant metric. Then,  from Theorem \ref{unica}, since any isometry must map the
canonical connection into itself, $\Iso(M) \subset \Aff (M, \nabla^c)$. But always $\Aff(M, \nabla^c) \subset
\Aff(M, \nabla)$ (affine transformations with respect to the Levi-Civita connection). In fact, any $g \in
\Aff(M, \nabla^c)$ maps Riemannian geodesics, which are the same as the $\nabla^c$-geodesics, into Riemannian
geodesics. Then, since $\nabla$ is torsion free, $g \in \Aff(M, \nabla)$, as it is well known (see \cite{R}).
Then
$$\Iso(M) \subset \Aff(M, \nabla^c) \subset \Aff(M, \nabla).$$
But, since $M$ is compact, $\Iso_0(M) = \Aff_0(M, \nabla)$ (see Remark \ref{affinecompact}). Hence,
$$\Iso_0(M) = \Aff_0(M, \nabla^c).$$

If  $M$ is isometric to a compact Lie group with with a bi-invariant metric then $\Iso_0(M) \subset \Aff_0(M,
\nabla^c)$ (see  Remark \ref{pruebaaffinecompact}). But, as observed before, $\Aff(M, \nabla^c) \subset \Aff(M,
\nabla)$. By making use of Remark \ref{affinecompact} we have,  also in this case, that
$$\Iso_0(M) = \Aff_0(M,\nabla ^c).$$
  The  geodesic symmetry, as observed in Remark \ref{pruebaaffinecompact}, does not preserve the
  canonical connection
 of the Lie group (if different from the Levi-Civita
connection). This concludes the proof.
\end{proof}

\begin{proof}[Proof of Corollary \ref{Ziller}] Observe first that $G_0$ is semisimple. In fact, let $\mathcal {D}$ be the distribution on
$M$ given by the tangent spaces to the orbits of the  maximal abelian normal (connected) subgroup $A$ of $G_0$.
Such a distribution must be $G$-invariant and so, since the action is effective, $\mathcal {D}_q = T_qM$, for
all $q\in M$. Then $A$ acts transitively on $M$ and so $M$ is flat, a contradiction. Let us endow $M$ with a
normal homogeneous metric $g'$, with respect to the decomposition $\cg = \gh \oplus \gh^\perp$. Namely,  the
scalar product in $T_pM \simeq \gh^\perp$ is the restriction of $-B$, where $B$ is the Killing form of $G_0$.

Such a nice metric $g'$ must  also be $G$-invariant, since any element of $H$ preserves both $\gh$ and $B$.
Since $H$ acts irreducibly on the tangent space, one has that $g'$, up to scaling, coincides with the metric $g$
of $M$.
 Let $\operatorname {Tr}(\nabla^c)$ be the normal subgroup of
$\Aff_0(M, \nabla^c)$ which consists of the transvection with respect to the canonical connection $\nabla ^c$
associated to the reductive decomposition $\cg = \gh \oplus \gh^\perp$. Recall that a  transvection is an
$\nabla^c$-affine transformation that preserves any holonomy subbundle of the orthonormal frame bundle. It is
well known that the Lie algebra of $\operatorname {Tr}(\nabla^c)$ is $\mathfrak{tr}(\nabla^c) = [\gh^\perp ,
\gh^\perp] + \gh^\perp$, not a direct sum in general,  which implies that $\operatorname {Tr}(\nabla^c) \subset
G_0$. By Theorem \ref{fullisometry}, $\Iso _0 (M) = \Aff_0(M, \nabla^c)$. Then $\operatorname {Tr}(\nabla^c)$ is
a normal subgroup of $\Iso _0 (M)$.

We shall prove that these groups coincide. In fact, assume that
 $\operatorname {Tr}(\nabla^c)$ is
properly contained in $ \Iso_0 (M)$. Let  $\cg '$ be a complementary ideal, in the Lie algebra of $\Iso _0(M)$,
of the ideal $\mathfrak{tr}(\nabla^c)$ (see \cite{R}). Then, if $0\neq X\in \cg '$, the field $z(q) = X . q$ is
$\operatorname {Tr}(\nabla^c)$-invariant. Then the isotropy of $\operatorname {Tr}(\nabla^c)$ fixes the vector
$z(p)$, where $p= [e]$ . But the subspace $\bbw$  of $T_pM$ which consist of the fixed vector of $\operatorname
{Tr}(\nabla^c)_p$ is invariant under $H$, since  $\operatorname {Tr}(\nabla^c)$ is a normal subgroup of $G$.
Recall that, from part (ii) of Theorem \ref{fullisometry}, since $M$ is not isometric to a compact Lie group,
one has that  $G\subset \Aff(M, \nabla^c)$. Since $H$ acts irreducibly, $\bbw = T_pM$. Then $\operatorname
{Tr}(\nabla^c)$ acts simply transitively on $M$. Then $[\gh^\perp , \gh^\perp] \subset \gh^\perp$ and therefore
$\gh^\perp$ is an ideal of $\cg$. This implies that $M$ is  isometric to the Lie group $\operatorname
{Tr}(\nabla^c)$ with a bi-invariant metric, a contradiction. Then $\operatorname {Tr}(\nabla^c) =  \Iso_0 (M)$
which implies that  $G_0 = \Iso _0 (M)$, since $\operatorname {Tr}(\nabla^c) \subset G_0 \subset \Iso _0 (M)$.
\end{proof}

\begin{nota} \label{strongly} For strongly isotropy irreducible spaces one needs not to assume that $M$
is not isometric to a compact Lie group with a bi-invariant metric. The proof is the same,  since always $\Iso
_0(M) = \Aff_0(M, \nabla^c)$, as it follows from Theorem \ref{fullisometry} (i) and we do not need to use part
(ii) of this result.
\end{nota}

\begin{nota} \label{affinecompact} (see  \cite{R}). Let $M$ be a compact Riemannian manifold and let $X$ be
an affine Killing field, i.e.\, the flow $\varphi _s$ associated to $X$ preserves the Levi-Civita connection
$\nabla$. Let $\gamma (t)$ be an arbitrary geodesic in $M$. Then $X(\gamma (t))$ is a Jacobi field along $\gamma
(t)$. In fact, for any $s\in \bbr$, $\gamma _s(t) = \varphi _s (\gamma (t))$ is a geodesic, since an affine
transformation maps geodesics into geodesics (note  that  $X (\gamma(t))  = \frac {\ \partial}{\partial
s}_{|0}\gamma _s(t)$). Since $M$ is compact, $X$ is bounded, and so, the projection of $ X (\gamma(t))$ to
$\gamma ' (t)$ is constant (since, it is of the form $(a+bt)\gamma ' (t)$). So, differentiating $\langle X
(\gamma(t)) , \gamma ' (t)\rangle$, one has that
$$\langle \nabla _ {\gamma ' (t)}  X (\gamma(t)) , \gamma ' (t)\rangle$$
which is the (Riemannian) Killing equation, since $\gamma $ is arbitrary. So, $X$ is a (Riemannian) Killing
field. This implies that $\Aff_0 (M, \nabla) = \Iso _0(M)$.
\end{nota}

\begin{nota} \label{projective}
Let $M^n =  G/G_p$ be a compact Riemannian manifold such that $(G_p)_0 = \SO (T_pM)$. Then $M$ is isometric to
the sphere $S^n$ or to the real projective space $\bbr P^n$. In fact, since the isotropy acts transitively on
the sets of planes of the tangent space at $p$, $M$ has constant sectional curvatures $\kappa$ ($\kappa >0$,
since $M$ is homogeneous and compact). Let  $S^n = \SO (n+1)/ \SO (n)$ be the universal cover of $M$ and let
$\Gamma$ be the deck transformations. Observe that the  isotropy $g\SO (n)g^{-1}$ of the sphere at $g . e_1$,
projects down to $M$. So, since $\Gamma$ is discrete  and the isotropies are connected, the deck transformation
group commutes with any isotropy of the sphere. This implies that $\Gamma$ commutes with $\SO(n+1)$. Therefore
$\Gamma = \{\Id\}$ or $\Gamma =\{ \Id ,  -\Id \}$.
\end{nota}

\begin{nota}
Let $M = G/G_p$ be a naturally reductive space. Let $\tilde M = \tilde G / \tilde G_{\tilde p}$ be its universal
cover, where $\tilde G$ is the (connected) lift of $G$ to $\tilde M$ and $\tilde p$ projects down to $p$.
Observe that $\tilde M$ is also a naturally reductive space. Since $\tilde G$ admits a bi-invariant metric, any
difeomorphism of $\tilde M$ that commutes with $\tilde G$ is volume preserving (see Remark \ref{volume}). Then,
with  the same arguments of the proof of   Theorem \ref{unica},  one has  that the canonical connection of
$\tilde M$ is unique, provided  $\tilde M$ does not split off either a sphere or a simple Lie group (with a
bi-invariant metric). So, in this case, the canonical connection of $M$ is also unique.
\end{nota}

\begin{nota}[Holonomy of naturally reductive spaces]

Let $M^n$ be a locally irreducible na\-tu\-ra\-lly reductive space with associated canonical connection
$\nabla^c$. The difference tensor $D := \nabla - \nabla^c$  between the Levi-Civita connection and the canonical
connection is totally skew. Moreover, this tensor gives the covariant derivative of Killing fields at a point $p
\in M$ (see the beginning of this section). Thus, from \cite{K}, \cite{AK},  for any  $v \in T_pM$, $D_v$
belongs to the (restricted) holonomy  algebra (see \cite{CDO}). If $M$ is not a symmetric space, then $D\neq 0$.
Moreover, by making use of  the strong skew-torsion holonomy theorem, the restricted holonomy group must be
$\SO(n)$. This result extends that  of  Wolf \cite{Wo} for strongly isotropy irreducible spaces.
\end{nota}


\section{The calculation of $\Aff(M, \nabla^c)$}

In this section we show how to compute geometrically the full $\nabla^c$-affine group (identity component) for a
naturally reductive space. We refer to \cite{R} for  more details. Let $M = G/H$ be a compact naturally
reductive space, with reductive decomposition $\cg = \gh \oplus \gm$ and associated canonical connection
$\nabla^c$. Assume that $H$ is the isotropy group at $p \in M$. We consider the transvection group of the
canonical connection, $\Tr(\nabla^c)$, which consists of all $\nabla^c$-affine transformations that preserves
the $\nabla ^c$-holonomy subbundles of the (orthogonal) frame  bundle. The transvection group is a normal
subgroup of $\Aff(M, \nabla^c)$. Moreover, $\Tr( \nabla^c)$ is contained in $G$ and the Lie algebra of $\Tr(
\nabla^c)$ is given by
$$\mathfrak{tr}(\nabla^c) = [\gm, \gm] + \gm$$
(not direct sum, in general).

Now, since $\Tr( \nabla^c)$ is a normal subgroup of $\Aff(M, \nabla^c)$, then $\mathfrak{tr}(\nabla^c)$ is an
ideal of $\mathfrak{aff}(M, \nabla^c)$, the Lie algebra of the affine group. Let us complement $\mathfrak{tr}(
\nabla^c)$ with another ideal (which corresponds to a normal subgroup of $\Aff(M, \nabla^c)$). That is, there
exists an ideal $\gg'$ of $\mathfrak{aff}(M, \nabla^c)$ such that
$$\mathfrak{aff}(M, \nabla^c) = \mathfrak{tr}( \nabla^c) \oplus \gg'.$$
Observe  that $\Aff_0(M, \nabla^c)$ is a closed, and hence compact, subgroup of $\Aff _0(M, \nabla) = \Iso _0
(M)$; see the proof of Theorem \ref{fullisometry} (last equality follows from Remark \ref{affinecompact}).

If $X \in \gg'$, then $\varphi_t$, the flow of $X$, commutes with $\Tr( \nabla^c)$ and therefore  $X$, regarded
as a field on $M$, is $\Tr( \nabla^c)$-invariant. Then $\gg'$ may be regarded as a subspace of
$\Tr(\nabla^c)$-invariant fields of $M$. Therefore, $\mathfrak{aff}(M, \nabla^c)$ is given by the transvection
algebra and a subspace of the $\Tr(\nabla^c)$-invariant fields. Now, we can write $\Aff_0(M, \nabla^c)$ as (not
direct) product of two groups of $\nabla ^c$-affine transformations. One can improve this presentation for a
quasi-direct product (i.e.\ discrete intersection of the factors, see \cite{R}).

\begin{nota}[see \cite{R}] If $M$ is a normal homogeneous space, one has that any $G$-invariant field lies in
$\mathfrak{aff}(M, \nabla^c)$. In fact, let $\varphi_t$ be the local flow of a $G$-invariant field $X$, then
$\varphi_t(p)$ is a fixed point of $H$, but the isotropy group does not change along the set of fixed points of
$H$. Thus,  the reductive decomposition in $\varphi_t(p)$ is the same, since $\gm = \gh^\bot$. From this one has
that $\varphi_t$ is $\nabla^c$-affine.
\end{nota}

\begin{nota}  It is  a well-known fact that a diffeomorphism is affine (i.e.\, it preserves some  connection) if
it maps geodesics into geodesics and it preserves the torsion tensor. For a canonical connection $\nabla ^c$ of
$M$, any linear isometry $\ell: T_pM \rightarrow  T_qM$  which maps the  canonical curvature and torsion tensors
tensors at $p$ into the same objects at $q$, extends to a $\nabla ^c$-affine transformation, which is also an
isometry. This is because the canonical curvature and torsion tensors are $\nabla ^c$-parallel (see the
introduction).
\end{nota}

\section{Appendix}

We shall give a geometric and classification free proof of the following result which follows from Joseph Wolf
classification of strongly isotropy irreducible spaces.

\begin{prop} [\cite{Wo}] \label{grouptype}  Let $M = G/H'$ be a Riemannian compact homogeneous manifold such that the connected
component $H$ of $H'$acts on $T_pM$, $p=[e]$, as the $\Ad$-representation of a compact simple Lie group. Then
$M$ is isometric to a compact simple Lie group with a bi-invariant metric.
\end{prop}

\begin{proof} If the rank of compact Lie group is $1$ then $M$ has dimension $3$ and $H' = \SO(T_pM)$
(see Remark \ref{rankonegroup}). Then, by Remark \ref{projective}, $M$ is isometric to the sphere $S^3$ or to
the projective space $\bbr P^3$, and the conclusion holds. So, let us assume that the rank is $k\geq 2$ (which
coincides with the codimension of the principal $H$-orbits).

Identifying $\gh = \Lie(H) \simeq T_pM$, we have that $H$, via the isotropy representation, acts as the
$\Ad$-action of $H$ on $\gh$. Let $0 \neq v \in \gh$. The normal space to the orbit $H . v$ is given by
$$ \nu_v (H . v) = \cc (v):= \{\xi \in \gh : [\xi,v]=0\}.$$
We have that $v = \Exp (tv) . v = \Ad (\Exp (tv))(v)$.  Then, for all $t \in \bbr$,    $\Ad (\Exp (tv))$ leaves
$\nu_v (H . v)$ invariant. Moreover, from the above equality, the set of fixed points of the one-parameter group
of linear isometries $\{\Ad (\Exp (tv))\}$ of  $\gh$ is just exactly the normal space $\nu_v (H . v)$ (see
Remark \ref{normal}).

Let us now write
$$dh_{t\vert _p} = \operatorname {Ad}(\Exp (tv))$$
where $h_t \in H$. Then $S := \{h_t\}$ is a one-parameter Lie group of isometries such that $M^v : =
\operatorname {exp}_p (\nu_v (H.v))$ is the connected component (containing $p$) of the set of fixed points of
$S$. Then $M^v$ is a totally geodesic submanifold of $M$. Observe that $M^v$, by Lemma \ref{fixedset},  is a
homogeneous submanifold of $M$. Moreover, as it is not hard to see, the isotropy algebra of $M^v$ is $\cc (v)$.
In the case that  $H.v$ is a maximal focal orbit,
$$\cc (v) = \mathbb{R}v\oplus \tilde \gh _v$$
where $\tilde \gh _v$ is a semisimple Lie subalgebra of $\gh$. The Lie algebra $\tilde \gh _v$ may be regarded
as the normal holonomy algebra of the maximal focal orbit $H . v$. The Lie algebra $\tilde \gh _v$ coincides
with its own normalizer in $\mathfrak {so}(\{v\}^\perp)$, since it acts as an $s$-representation, where
$\{v\}^\perp$ is regarded inside  $\nu _v (H . v)$ (see \cite{BCO}, p.\ 192). From this property, since $\tilde
\gh _v$ is semisimple, one has that
$$\{ 0 \} = \cc(\tilde \gh _v) = \{x \in \mathfrak{so}(\{v\}^\perp): [x, \tilde \gh _v] = 0\}.$$

So, we are under the assumptions of Lemma \ref{3.3}. Then, $M^v$ splits, locally, the geodesic $\gamma _v$.

Let now $H . u$ be a principal orbit and choose $0\neq  v\in u + \nu _u(H . u)= \nu _u(H . u)$ such that $H . v$
is a most singular orbit (i.e.\ $v$ belongs to a one dimensional simplex of a Weyl chamber of the normal space
$\nu _u(H . u)$; recall that the principal $H$-orbits are isoparametric submanifolds \cite{PT, BCO}). Since
$M^v$ splits locally $\gamma _v$, then $M^u$ splits also locally $\gamma _v$, since $M^v$ is a totally geodesic
submanifold of $M^u$ (due to the fact that $\nu _u(H . u) \subset \nu _v(H . v)$).

Let $W$ be the Weyl group of $\nu_u(H . u)$ which acts irreducibly on $ u + \nu _u(H . u) = \nu _u(H . u)$.
Given $g \in W$ there exists $h\in H$ such that $h . \nu _u(H . u) = \nu _u(H . u)$ and $h_{\vert \nu _u(H . u)}
= g$. Then $M^u = h . M^u$ splits off  also $\gamma _{h . v} = \gamma _{g(v)}$, for all $g \in W$. Then $M^u$ is
flat (and compact), since $W.v$ generates $\nu _u(H . u)$. Let now $z \in T_pM$ be arbitrary. There exists $w
\in H.u$ such that $z \in \nu_w (H . u) = \nu_w (H . w)$ (for instance, by choosing  a point $w$ where the
height function $x\mapsto \langle x,z \rangle$, $x\in H . u$, achieves its maximum value). Then the arbitrary
geodesic $\gamma _z$, by $p$, is contained in the compact $k$-flat $M^w$, where $k = \text {rank}(H)\geq 2$.
Since $M$ is homogeneous, we conclude that any geodesic in $M$ must be contained in a compact flat. Then, by
\cite{HPTT}, $M$ is globally symmetric of rank at least two (see \cite{EO} for a conceptual proof).

Observe that $H$, regarded as a subgroup of the isotropy of $M$,  coincides with  the connected component of the
full isotropy $\Iso_0(M)_p$ of $M$ at $p$. If not, $\Iso_0(M)_p$ would be transitive on the sphere, by Simons'
Holonomy Theorem \cite{S,O2} (since $s$-representations can not be enlarged without acting transitively on the
sphere). But then the symmetric space $M$ would be of rank one,  a contradiction since there exists non trivial
flats in $M$.

So, since $M$ is symmetric, $H = \Iso_0(M)_p  $ coincides with the (restricted) holonomy  group of $M$ at $p$
(via the isotropy representation).

Let now $X$ be the symmetric space $H$ endowed with a bi-invariant metric. Identifying $T_pM \simeq T_eX$ one
has both $M$ and $X$ have the same (restricted) holonomy. Then,  the curvature tensors of $M$ and $X$ differ by
a positive scalar multiple (this follows for instance from Lemma 3.3 in \cite{O2}).

Thus $M$ is locally a symmetric space of group type. Since $M$ is globally symmetric it is not hard to see  that
$M$ must be globally isometric to a (connected and compact) Lie group with a bi-invariant metric.
\end{proof}

\end{document}